\documentstyle[twoside,amstex]{article}
\input amssym.def
\input amssym.tex
\def\bc{\begin{center}}
\def\ec{\end{center}}
\def\no{\noindent}
\def\hang{\hangindent\parindent}
\def\textindent#1{\indent\llap{[#1]\enspace}\ignorespaces}
\def\re{\par\hang\textindent}
\begin{document}
\thispagestyle{empty} \vspace*{3 true cm} \pagestyle{myheadings}
\markboth {\hfill {\sl Huanyin Chen and Marjan Sheibani}\hfill}
{\hfill{\sl Elementary Matrix Reduction Over B\'{e}zout Duo-domains}\hfill} \vspace*{-1.5 true cm} \bc{\large\bf Elementary Matrix Reduction Over B\'{e}zout Duo-domains}\ec

\vskip6mm
\bc{{\bf Huanyin Chen}\\[2mm]
Department of Mathematics, Hangzhou Normal University\\
Hangzhou 310036, China\\
huanyinchen@@aliyun.com}\ec

\bc{{\bf Marjan
Sheibani}\\[2mm]
Faculty of Mathematics, Statistics and Computer Science\\
Semnan University, Semnan, Iran\\
m.sheibani1@@gmail.com}\ec

\vskip10mm
\begin{minipage}{120mm}
\no {\bf Abstract:} A ring $R$ is an elementary divisor ring if every matrix over $R$ admits a diagonal reduction. If $R$ is
an elementary divisor domain, we prove that $R$ is a B\'{e}zout duo-domain if and only if
for any $a\in R$, $RaR=R\Longrightarrow ~\exists ~s,t\in R$ such that $sat=1$. We further explore various stable-like conditions on a B\'{e}zout duo-domain under which it is an elementary divisor domain. Many known results are thereby generalized to much wider class of rings, e.g. [3, Theorem 3.4.], [5, Theorem 14], [9, Theorem 3.7], [13, Theorem 4.7.1] and [14, Theorem 3].
\vskip3mm \no {\bf MR(2010) Subject Classification}: 15A21, 16S99, 16L99, 19B10.
\end{minipage}

\vskip15mm \bc{\bf 1. Introduction}\ec

\vskip4mm \no Throughout this paper, all rings are associative with an identity. A matrix $A$ (not necessarily square) over a ring
$R$ admits diagonal reduction if there exist invertible matrices
$P$ and $Q$ such that $PAQ$ is a diagonal matrix $(d_{ij})$, for
which $d_{ii}$ is a full divisor of $d_{(i+1)(i+1)}$ (i.e., $Rd_{(i+1)(i+1)}R\subseteq d_{ii}R\bigcap Rd_{ii}$) for each $i$. A
ring $R$ is called an elementary divisor ring provided that every
matrix over $R$ admits a diagonal reduction. A ring $R$ is B\'{e}zout ring if every finitely generated right (left) ideal is principal.
Clearly, every elementary divisor ring is a B\'{e}zout ring.
An attractive problem is to investigate various conditions under which a B\'{e}zout ring is an elementary divisor ring.

Commutative elementary divisor domains have been studies by many authors, e.g. [3],[5], [9-10], and [13-14]. But the structure of such rings in the noncommutative case has not been sufficiently studied (cf. [6] and [15]). A ring $R$ is duo if every right (left) ideal of $R$ is a two-sided ideal.
Obviously, every commutative ring is duo, but the converse is not true. In this paper, we are concern on when a B\'{e}zout duo-domain is an elementary divisor domain. Here, a ring $R$ is a domain ring if there is no any nonzero zero divisor.

In Section 2, We prove that an elementary divisor domain is a B\'{e}zout duo-domain if and only if
for any $a\in R$, $RaR=R\Longrightarrow ~\exists ~s,t\in R$ such that $sat=1$. A ring $R$ has stable range 1 if, $aR+bR=R$ with $a,b\in R\Longrightarrow \exists y\in R$ such that $a+by\in R$ is invertible. This condition plays an important rule in algebra K-theory. We refer the reader to [1] for the general theory of stable range 1. Further, we shall introduce certain stable-like conditions on B\'{e}zout duo-domains so that they are elementary divisor domains.
In Section 3, we prove that every locally stable B\'{e}zout duo-domain and every Zabavsky duo-domainis are
elementary
divisor domains (see Theorem 3.2 and Theorem 3.9). In Section 4,
we prove that every B\'{e}zout duo-domain of Gelfand range 1 is an
elementary divisor ring (see Theorem 4.3). Finally, in the last section, we establish related theorems on B\'{e}zout duo-domains by means of countability conditions on their maximal ideals. Many known results are
thereby generalized to much wider class of rings, e.g. [3, Theorem 3.4.], [5, Theorem 14], [9, Theorem 3.7], [13, Theorem 4.7.1] and [14, Theorem 3].

We shall use $J(R)$ and $U(R)$ to denote the Jacobson radical of $R$ and the set of all units in $R$, respectively. For any $a,b\in R$, $a|b$ means that $a\in bR$.

\vskip15mm\bc{\bf 2. B\'{e}zout Duo-domains}\ec

\vskip4mm A ring $R$ is right (left) quasi-duo if every right (left) maximal ideal of $R$ is an ideal. A ring $R$ is quasi-duo if it is both right and left quasi-duo. It is open whether there exists a right quasi-duo ring which is not left
quasi-duo [8]. We have

\vskip4mm \hspace{-1.8em} {\bf Theorem 2.1.}\ \ {\it Let $R$ be an elementary divisor ring. Then the following are equivalent:}
\vspace{-.5mm}
\begin{enumerate}
\item [(1)] {\it $R$ is quasi-duo.}
\vspace{-.5mm}
\item [(2)] {\it $R$ is left (right) quasi-duo ring.}
\vspace{-.5mm}
\item [(3)] {\it For any $a\in R$, $RaR=R\Longrightarrow ~\exists ~s,t\in R$ such that $sat=1$.}
\end{enumerate}
\vspace{-.5mm} {\it Proof.}\ \ $(1)\Rightarrow (2)$ This is trivial.

$(2)\Rightarrow (3)$ We may assume that $R$ is right quasi-duo. Suppose that $RaR=R$ with $a\in R$. Write $1=r_1as_1+\cdots +r_nas_n$. Then $r_1aR+\cdots +r_naR=R$. In light of [8, Theorem 3.2], $Rr_1a+\cdots +Rr_na=R$. This shows that $sa=1$ for some $s\in R$, as desired.

$(3)\Rightarrow (1)$ Suppose that $aR+bR=R$ with $a,b\in R$. Let $A=\left(
\begin{array}{cc}
a&0\\
b&0
\end{array}
\right)$. Since $R$ is an elementary divisor ring, there exist $U,V\in GL_2(R)$ such that
$$UAV=\left(
\begin{array}{cc}
w&0\\
0&v
\end{array}
\right),$$ where $RvR\subseteq wR\bigcap Rw$. It follows from $aR+bR=R$ that $RwR+RvR=R$, and so $RwR=R$.
Since $R$ satisfies Zabavsky condition, we can find some $s,t\in R$ such that
$swt=1$. Then we
verify that
$$
\left(
\begin{array}{cc}
s&0\\
1-wts&wt\end{array}
\right)\left(
\begin{array}{cc}
0&w\\
v&0
\end{array}
\right)\left(
\begin{array}{cc}
sw&0\\
1-tsw&t
\end{array}
\right)
=\left(
\begin{array}{cc}
*&1\\
**&*\end{array}
\right),$$
so we show that
$$
\left(
\begin{array}{cc}
s&0\\
1-wts&wt\end{array}
\right)UAV\left(
\begin{array}{cc}
0&1\\
1&0
\end{array}
\right)
\left(
\begin{array}{cc}
sw&0\\
1-tsw&t
\end{array}
\right)\left(
\begin{array}{cc}
0&1\\
1&0
\end{array}
\right)
=\left(
\begin{array}{cc}
1&*\\
**&*\end{array}
\right).$$
One easily checks that
$$\begin{array}{c}\left(
\begin{array}{cc}
s&0\\
1-wts&wt\\
\end{array}
\right)=B_{21}(-wt)B_{12}(s-1)B_{21}(1)B_{12}(wt-1);\\
\left(
\begin{array}{cc}
sw&0\\
1-tsw&t\\
\end{array}
\right)=B_{21}(-t)B_{12}(sw-1)B_{21}(1)B_{12}(t-1).
\end{array}
$$ Here, $$B_{12}(*)=\left(
\begin{array}{cc}
1&*\\
0&1\\
\end{array}
\right), B_{21}(*)=\left(
\begin{array}{cc}
1&0\\
**&1\\
\end{array}
\right).$$
Thus, $$
\left(
\begin{array}{cc}
s&0\\
1-wts&wt\\
\end{array}
\right), \left(
\begin{array}{cc}
sw&0\\
1-tsw&t\\
\end{array}
\right) \in GL_2(R).$$ Therefore we can find some $P=(p_{ij}),Q=(q_{ij})\in GL_2(R)$ such that
$$PAQ=\left(
\begin{array}{cc}
1&0\\
0&u
\end{array}
\right).$$ This implies that $(p_{11}a+p_{12}b)q_{11}=1$, and then $Ra+Rb=R$. It follows by [8, Theorem 3.2] that $R$ is left quasi-quo ring. Likewise, $R$ is right quasi-duo ring.\hfill$\Box$

\vskip4mm Following Zabavsky [13], we say that a ring $R$ satisfies Lam condition if $RaR=R\Longrightarrow ~a\in R$ is invertible.

\vskip4mm \hspace{-1.8em} {\bf Corollary 2.2.}\ \ {\it Let $R$ be an elementary divisor ring. Then the following are equivalent:}
\vspace{-.5mm}
\begin{enumerate}
\item [(1)] {\it $R$ is quasi-duo.}
\vspace{-.5mm}
\item [(2)] {\it $R$ satisfies Lam condition.}
\end{enumerate}
\vspace{-.5mm} {\it Proof.}\ \ $(1)\Rightarrow (2)$ This is obvious.

$(2)$ If $RaR=R$, by hypothesis, $a\in U(R)$, which completes the proof by Theorem 2.1.\hfill$\Box$

\vskip4mm Recall that a ring $R$ satisfies Dubrovin condition if, for any $a\in R$ there exists $b\in R$ such that $RaR=bR=Rb$. We now extend [13, Theorem 4.7.1] as follows:

\vskip4mm \hspace{-1.8em} {\bf Corollary 2.3.}\ \ {\it Let $R$ be an elementary divisor ring. If $RaR=R\Longrightarrow ~\exists ~s,t\in R$ such that $sat=1$, then $R$ satisfies Dubrovin condition.}
\vskip2mm\hspace{-1.8em} {\it Proof.}\ \ Let $0\neq a\in R$. Since $R$ is an elementary divisor ring, there exist $U=(u_{ij}), V=(v_{ij})\in GL_2(R)$ such that
$$\left(
\begin{array}{cc}
a&0\\
0&a
\end{array}
\right)U=V\left(
\begin{array}{cc}
b&0\\
0&c
\end{array}
\right),$$ where $RcR\subseteq bR\bigcap Rb.$ It follows that $au_{12}=v_{12}c$ and $au_{22}=v_{22}c$. Since $U\in GL_2(R)$, we see that
$Ru_{12}+Ru_{22}=R$. In view of Theorem 2.1, $R$ is quasi-duo. By virtue of [8, Theorem 3.2], we have $u_{12}R+u_{22}R=R$. Write
$$u_{12}x+u_{22}y=1~\mbox{for some}~x,y\in R.$$ Then $$a=au_{12}x+au_{22}y=v_{12}cx+v_{22}cy\in RcR;$$ hence, $$RaR\subseteq RcR.$$ But $$RcR\subseteq bR\bigcap Rb\subseteq RbR\subseteq RaR.$$ Thus, $$RaR=RcR.$$ Clearly, $$RcR\subseteq bR\subseteq RbR\subseteq RaR=RcR~\mbox{and}~RcR\subseteq Rb\subseteq RbR\subseteq RaR=RcR.$$ Therefore $RaR=RcR=bR=Rb$, as desired.\hfill$\Box$

\vskip4mm \hspace{-1.8em} {\bf Corollary 2.4 [13, Theorem 4.7.1].}\ \ {\it Let $R$ be an elementary divisor ring. If $R$ satisfies Lam condition, then $R$ satisfies Dubrovin condition.}
\vskip2mm\hspace{-1.8em} {\it Proof.}\ \ Since $R$ satisfies Lam condition, we see that $RaR=R$ implies $a\in R$ is invertible. In view of Corollary 2.3,
 $R$ satisfies Dubrovin condition, as asserted.\hfill$\Box$

\vskip4mm \hspace{-1.8em} {\bf Theorem 2.5.}\ \ {\it Let $R$ be an elementary divisor domain. Then the following are equivalent:}
\vspace{-.5mm}
\begin{enumerate}
\item [(1)] {\it $R$ is a B\'{e}zout duo-domain.}
\vspace{-.5mm}
\item [(2)] {\it $R$ satisfies Lam condition.}
\vspace{-.5mm}
\item [(3)] {\it For any $a\in R$, $RaR=R\Longrightarrow ~\exists ~s,t\in R$ such that $sat=1$.}
\end{enumerate}
\vspace{-.5mm} {\it Proof.}\ \ $(1)\Rightarrow (2)\Rightarrow (3)$ These are obvious.

 $(3)\Rightarrow (1)$ In view of Theorem 2.1, $R$ is quasi-duo. In view of Corollary 2.3, $R$ satisfies Dubrovin condition.
 Therefore $R$ is a duo ring, in terms of [15, Theorem 1].\hfill$\Box$

\vskip4mm An element $e\in R$ is infinite if there exist orthogonal
idempotents $f,g\in R$ such that
$e=f+g$ while $eR\cong fR$ and $g\neq 0$. A simple domain is said to
be purely infinite if every nonzero
left ideal of $R$ contains an infinite idempotent, As is well known, a ring $R$ is a purely infinite simple ring if and only if it is not a division ring and for any nonzero $a\in R$ there exist $s,t\in R$ such that $sat=1$ (cf. [1]).

\vskip4mm \hspace{-1.8em} {\bf Corollary 2.6.}\ \ {\it Every purely infinite simple domain is not a B\'{e}zout domain.}
\vskip2mm\hspace{-1.8em} {\it Proof.}\ \ Let $R$ be a purely infinite ring. Suppose that $R$ is a B\'{e}zout domain. In view of [13, Theorem 1.2.6], $R$ is an Hermite ring.
Thus, every $1\times 2$ and $2\times 1$ matrices over $R$ admit a diagonal reduction. Let $A=(a_{ij})\in M_2(R)$. There exists $Q\in GL_2(R)$ such that
$$AQ=\left(
\begin{array}{cc}
a&0\\
b&c
\end{array}
\right)$$ for some $a,b,c\in R$. If $a=b=c=0$, then $AQ=0$ is a diagonal matrix. Otherwise, we may assume that $b\neq 0$. Thus, we can find $s,t\in R$ such that
$sbt=1$. Then we
check that
$$\left(
\begin{array}{cc}
sb&0\\
1-tsb&t
\end{array}
\right)\left(
\begin{array}{cc}
a&0\\
b&c
\end{array}
\right)
\left(
\begin{array}{cc}
s&1-bts\\
0&bt\end{array}
\right)
=\left(
\begin{array}{cc}
*&*\\
1&*\end{array}
\right),$$
so we have $$B_{21}(*)\left(
\begin{array}{cc}
0&1\\
1&0\end{array}
\right)\left(
\begin{array}{cc}
sb&1-tsb\\
0&t
\end{array}
\right)\left(
\begin{array}{cc}
a&0\\
b&c
\end{array}
\right)
\left(
\begin{array}{cc}
s&0\\
1-bts&bt\end{array}
\right)B_{21}(*)=\left(
\begin{array}{cc}
1&0\\
0&*\end{array}
\right).$$
As in the proof of Theorem 2.1, we see that
$$\left(
\begin{array}{cc}
s&0\\
1-bts&bt\\
\end{array}
\right), \left(
\begin{array}{cc}
sb&1-tsb\\
0&t\\
\end{array}
\right)\in GL_2(R).$$ Therefore, $A$ is equivalent to a diagonal matrix $diag(1,*)$. Hence, $R$ is an elementary divisor ring.
If $RaR=R$, then $a\neq 0$, and so $sat=1$ for some $s,t\in R$. In light of Theorem 2.1, $R$ is a quasi-duo ring. Thus, $R$ is Dedekind-finite, i.e., $uv=1$ in $R\Longrightarrow vu=1$. Let $0\neq x\in R$. Then there exists $c,d\in R$ such that $cxd=1$. Hence, $xdc=dcx=1$. This implies that $x\in U(R)$, and thus $R$ is a division ring, a contradiction. Therefore we conclude that $R$ is not a B\'{e}zout domain.\hfill$\Box$

\vskip4mm \hspace{-1.8em} {\bf Lemma 2.7.}\ \ {\it Let $R$ be a B\'{e}zout duo-domain. Then the following are equivalent:}
\vspace{-.5mm}
\begin{enumerate}
\item [(1)] {\it $R$ is an elementary divisor ring.}
\vspace{-.5mm}
\item [(2)] {\it $aR+bR+cR=R$ with $a,b,c\in R\Longrightarrow ~\exists p,q\in R$ such that $(pa+qb)R+qcR=R$.}
\end{enumerate}
\vspace{-.5mm} {\it Proof.}\ \ This is clear, by [15, Lemma 2 and Lemma 3].\hfill$\Box$

\vskip4mm \hspace{-1.8em} {\bf Theorem 2.8.}\ \ {\it Let $R$ be a B\'{e}zout duo-domain. Then the following are equivalent:}
\vspace{-.5mm}
\begin{enumerate}
\item [(1)] {\it $R$ is an elementary divisor ring.}
\vspace{-.5mm}
\item [(2)] {\it $Ra+Rb+Rc=R$ with $a,b,c\in R\Longrightarrow \exists r,s,t\in R$ such that $s|rt$ and $ra+sb+tc=1$.}
\end{enumerate}\vspace{-.5mm} {\it Proof.}\ \ $(1)\Rightarrow (2)$ Suppose $Ra+Rb+Rc=R$ with $a,b,c\in R$. Then $aR+bR+cR=R$, as $R$ is duo. In view of Lemma 2.7, there exist $p,q\in R$ such that $(pa+qb)R+qcR=R$. Hence, $R(pa+qb)+R(qc)=R$. Write $1=x(pa+qb)+y(qc)$. Then $(xp)a+(xq)b+(yq)c=1$. As $R$ is a duo-domain, we see that $(xp)(yq)\in xRq=Rxq$, as desired.

$(2)\Rightarrow (1)$ Let $a,b,c\in R$ be such that $aR+bR+cR=R$. Then $Ra+Rb+Rc=R$. By hypothesis, there exist $r,s,t\in R$ such that $tr\in Rs$ and $ta+sb+rc=1$.
Since $R$ is a B$\acute{e}$zout ring, we can find a $d\in R$ such that $Rd=dR=Rs+Rt=sR+tR$, as $R$ is a duo-domain.

Case I. $d=0$. Then $c\in U(R)$, and so $(a+c^{-1}\cdot b)R+(c^{-1}\cdot c)R=R$.

Case II. $d\neq 0$. Write $s=dq$ and $t=dp$ for some $q,p\in R$. Then $Rdr=Rsr+Rtr\subseteq Rs=sR$, and so $dRr=Rdr\supseteq sR=dqR$.
This shows that $r\in qR=Rq$.
Therefore $1=ta+sb+rc=d(pa+qb)+kqc\in R(pa+qb)+R(qc)=(pa+qb)R+(qc)R$; hence, $(pa+qb)R+qcR=R$.
In light of Lemma 2.7, $R$ is an elementary divisor ring.\hfill$\Box$

\vskip4mm \hspace{-1.8em} {\bf Example 2.9.}\ \ Let $R$ be the ring of all quaternions. Then $R$ is a noncommutative division ring. Thus, $R$ is an elementary divisor duo-domain, but it is not commutative.

\vskip4mm \hspace{-1.8em} {\bf Theorem 2.10.}\ \ {\it Let $R$ be a B\'{e}zout duo-domain. Then $R$
is an elementary divisor ring if and only if $R/J(R)$ is an elementary divisor ring.}
\vskip2mm\hspace{-1.8em} {\it Proof.}\ \ $\Longrightarrow$ This is obvious.

$\Longleftarrow$ Suppose that $aR+bR+cR=R$ with $a,b,c\in R$. Then $\left(
\begin{array}{cc}
\overline{a}&\overline{0}\\
\overline{b}&\overline{c}
\end{array}
\right)\in M_2(R/J(R))$ admits a diagonal reduction. Thus, we can find some $(\overline{p_{ij}}), (\overline{q_{ij}})\in GL_2(R/J(R))$ such that
$$(\overline{p_{ij}})\left(
\begin{array}{cc}
\overline{a}&\overline{0}\\
\overline{b}&\overline{c}
\end{array}
\right)(\overline{q_{ij}})=diag(\overline{u}, \overline{v}),$$ where $\overline{RvR}\subseteq \overline{uR}\bigcap \overline{Ru}.$
As $R$ is a duo-domain, we see that $\overline{uR+vR}=\overline{R}$. Hence, $\overline{u}\in R/J(R)$ is right invertible. This implies
that $u\in R$ is invertible, as $R$ is quasi-duo. Therefore $(p_{11}a+p_{12}b)q_{11}+p_{12}cq_{12}=u$, and so $(p_{11}a+p_{12}b)R+p_{12}cR=R.$
In light of Lemma 2.7, $R$ is an elementary divisor domain.
\hfill$\Box$

\vskip4mm \hspace{-1.8em} {\bf Corollary 2.11.}\ \ {\it Let $R$ be a B\'{e}zout duo-domain. Then $R$
is an elementary divisor ring if and only if $R[[x_1,\cdots,x_n]]$ is an elementary divisor ring.}
\vskip2mm\hspace{-1.8em} {\it Proof.}\ \ $\Longleftarrow$ This is obvious, as $R$ is a homomorphic image of $R[[x_1,\cdots,x_n]]$.

$\Longrightarrow$ Let $\phi: R[[x_1,\cdots,x_n]]\rightarrow R$ is defined by $\phi(f(x_1,\cdots,x_n))=f(0,\cdots,0)$. It is obvious that $I:= Ker(\phi)\subseteq J(R[[x_1,\cdots,x_n]])$ and $R[[x_1,\cdots,x_n]]/I \cong R$. Since $R/I/J(R)/I\cong R/J(R)$, $R/I/J(R)/I$ is an elementary divisor ring if and only if so is $R/J(R)$. Since $R/I\cong R$, the result follows by Theorem 2.10.
\hfill$\Box$

\vskip15mm\bc{\bf 3. Locally Stable domains}\ec

\vskip4mm We say that an element $a$ of a duo-ring $R$ is stable if, $R/aR$ has stable range 1. A duo ring $R$ is locally stable if, $aR+bR=R$ with $a,b\in R\Longrightarrow ~\exists~ y\in R$ such that $a+by\in R$ is stable. For instance, every duo-ring of stable range 1 is locally stable. The purpose of this section is to investigate matrix diagonal reduction over locally stable B\'{e}zout duo-domain.

\vskip4mm \hspace{-1.8em} {\bf Lemma 3.1.}\ \
{\it Let $R$ be a B\'{e}zout domain. If $pR+qR=R$ with $p,q\in R$, then there exists a matrix
$\left(\begin{array}{cc}
p&q\\
**&*
\end{array}
\right)\in GL_2(R)$.}
\vskip2mm\hspace{-1.8em} {\it Proof.}\ \ As every B\'{e}zout domain is a Hermite domain, then there exists some $Q\in GL_2(R)$ such that $(p,q)Q=(d,0)$ for some $d\in R$. We have $p,q\in dR$ and also $R=pR+qR\subseteq dR$ this implies that $d$ is right invertible so there exists some $x\in R$ such that $dx=1$, now $xdx=x$ then $(xd-1)x=0$. Since $R$ is a domain, then $xd-1=0$ that implies $xd=1$ and so $d$ is invertible. Now $(p,q)Q\left(
\begin{array}{cc}
d^{-1}&0\\
0&d^{-1}
\end{array}
\right)=(1,0)$. Set $U^{-1}=Q\left(
\begin{array}{cc}
d^{-1}&0\\
0&d^{-1}
\end{array}
\right)$. Then $(p,q)=(1,0)U$, this means that $(p,q)$ is the first row of $U$ So $U=\left(
\begin{array}{cc}
p&q\\
**&*
\end{array}
\right)\in GL_2(R)$ as required. \hfill$\Box$

\vskip4mm Let $R$ be a B\'{e}zout domain. If $Rx+Ry=R$ with $x,y\in R$, similarly, we can find a matrix
$\left(\begin{array}{cc}
x&*\\
y&*
\end{array}
\right)\in GL_2(R)$.

\vskip4mm \hspace{-1.8em} {\bf Theorem 3.2.}\ \ {\it  Every locally stable B\'{e}zout duo-domain is an elementary divisor ring.}\vskip2mm\hspace{-1.8em} {\it Proof.}\ \
Let $R$ be a locally stable B\'{e}zout duo-domain.
Suppose that $aR+bR+cR=R$ with $a,b,c\in R$. Then $ax+by+cz=1$ for some $x,y,z\in R$.
Since $R$ is locally stable, there exists some $s\in R$ such that $R/(b+(ax+cz)s)R$ has stable range 1. Set $w=b+(ax+cz)s$. Then $\overline{aR+cR}=\overline{1}$ in $R/wR$. Thus, we have
$$\left(
\begin{array}{cc}
0&1\\
1&0
\end{array}
\right) \left(
\begin{array}{cc}
1&0\\
xz&1
\end{array}
\right) \left(
\begin{array}{cc}
a&0\\
b&c
\end{array}
\right)\left(
\begin{array}{cc}
1&0\\
zs&1
\end{array}
\right)=\left(
\begin{array}{cc}
w&c\\
a&0
\end{array}
\right).$$ In light of [13, Theorem 1.2.6], $R$ is a Hermite ring. Thus, we have a $Q\in GL_2(R)$ such that $(w,c)Q=(v,0)$, and then
$$\left(
\begin{array}{cc}
w&c\\
a&0
\end{array}
\right)W=\left(
\begin{array}{cc}
v&0\\
k&l
\end{array}
\right).$$ Clearly, $w\in vR$; hence, $R/vR\cong \frac{R/wR}{vR/wR}$ has stable range 1.

One easily checks that $vR+kR+lR=R$, and so $\overline{kR+lR}=\overline{R}$. Thus, we have a $t\in R$ such that $vR+(k+lt)R=R$.
Since $R$ is duo, we see that $Rv+(k+lt)R=R$. Write $pv+p'(k+lt)=1$. Then $$\left(
\begin{array}{cc}
v&0\\
k&l
\end{array}
\right)\left(
\begin{array}{cc}
1&0\\
t&1
\end{array}
\right)=\left(
\begin{array}{cc}
v&0\\
k+lt&1
\end{array}
\right).$$

By virtue of Lemma 3.1, we can find some $P,Q\in GL_2(R)$ such that $$P\left(
\begin{array}{cc}
a&0\\
b&c
\end{array}
\right)Q=\left(
\begin{array}{cc}
1&0\\
0&d
\end{array}
\right).$$ Write $P=\left(
\begin{array}{cc}
p&q\\
p'&q'
\end{array}
\right)$ and $Q^{-1}=\left(
\begin{array}{cc}
m&n\\
m'&n'
\end{array}
\right)$. Then $pa+qb=m$ and $qc=n$. As $mR+nR=R$, we see that $(pa+qb)R+qcR=R$. In light of Lemma 2.7, we complete the proof.\hfill$\Box$

\vskip4mm As an immediate consequence, we extend [9, Theorem 3.7] from commutative case to duo rings.

\vskip4mm \hspace{-1.8em} {\bf Corollary 3.3.}\ \ {\it Let $R$ be a B\'{e}zout duo-domain. If $R/aR$ has stable range 1 for all nonzero $a\in R$, then $R$ is an elementary divisor domain.}
\vskip2mm\hspace{-1.8em} {\it Proof.}\ \ Given $aR+bR=R$ with $a,b\in R$, then $a\neq 0$ or $a+b\neq 0$. Hence, $a+by\neq 0$ for some $y\in R$. By hypothesis, $\frac{R}{(a+by)R}$ has stable ring 1, i.e., $R$ is locally stable. Therefore $R$ is an elementary divisor domain, by Theorem 3.2.\hfill$\Box$

\vskip4mm Since every B\'{e}zout duo-domain of stable range 1 satisfies the condition in Corollary 3.3, we derive

\vskip4mm \hspace{-1.8em} {\bf Lemma 3.4.}\ \ Every B\'{e}zout duo-domain of stable range 1 is an elementary divisor domain.

\vskip4mm \hspace{-1.8em} {\bf Theorem 3.5.}\ \ {\it Let $R$ be a B\'{e}zout duo-domain. If $aR+bR=R\Longrightarrow \exists ~x,y\in R$ such that $xR+yR=R, xy\in J(R), a|x,b|y$, then $R$
is an elementary divisor ring.}
\vskip2mm\hspace{-1.8em} {\it Proof.}\ \ Suppose that $aR+bR=R$ with $a,b\in R$. Then there exist $x,y\in R$ such that $xR+yR=R, xy\in J(R), a|x,b|y$. Write
$xr+ys=1$ for some $r,s\in R$. Since $R$ is a duo-domain, we see that $xr\in Rx$, and then $(xr)^2\equiv xr (mod~J(R))$.

Set $g=xr$ and
$e=g+gxr(1-g)$. Then $e\in xR\subseteq aR$. One easily
checks that $1-e=\big(1-gxr\big)(1-g)\equiv 1-g (mod~J(R))$. Thus, we have some $d\in J(R)$ such that $1-e-d=1-g=ys\in yR\subseteq bR$.
Set $f=1-e-d$. Then $e+f=1-d\in U(R)$, and so $e(1-d)^{-1}+f(1-d)^{-1}=1$. This implies that $eR+fR=R$. One easily checks that
$a|e$ and $b|f$. Moreover $ef=(e-e^2)-ed\in J(R)$. Therefore $R/J(R)$ is a duo exchange ring. In view of [1, Corollary 1.3.5], $R/J(R)$ has stable range 1, and then so does $R$. According to Lemma 3.4, we complete the proof.\hfill$\Box$

\vskip4mm Recall that a ring $R$ is feckly clean if for any $a\in R$, there exists $e\in R$ such that $a-e\in U(R)$ and $eR(1-e)\subseteq J(R)$ [2]. We now derive

\vskip4mm \hspace{-1.8em} {\bf Corollary 3.6.}\ \ {\it Every feckly clean B\'{e}zout duo-domain is an elementary divisor ring.}
\vskip2mm\hspace{-1.8em} {\it Proof.}\ \ Let $R$ be a feckly clean B\'{e}zout duo-ring. Suppose that $aR+bR=R$ with $a,b\in R$. Then $ax+by=1$ for some $x,y\in R$.
Set $c=ax$. Then there exist $f\in R$ and
$u\in U(R)$ such that $c=f+u$ and $fR(1-f)\in J(R)$. Let
$g=(1-f)+(f-f^2)u^{-1}$. Then $g-g^2\equiv f-f^2\equiv 0
\big(~mod~J(R)\big)$. Hence, $g-g^2\in J(R)$. As $fu(1-f)\in J(R)$, we see that $fu-fuf\in J(R)$.
On the other hand, $(R(1-f)ufR)^2=R(1-f)u(fR(1-f))ufR\subseteq J(R)$, and then $(1-f)uf\in J(R)$. This implies that $uf=fuf$. Therefore,
$fu-uf\in J(R)$. Write $uf=fu+r$ for some $r\in J(R)$. Then there exists some $r'\in J(R)$ such that $a-g=f+u-1+f-(f-f^2)u^{-1}=(u-2fu-u^2+f-f^2)(-u^{-1})=(c-c^2)(-u^{-1})+r'$. Write $e=g-r'$. Then $a-e\in
(c-c^2)R$. Write $c-e=(c-c^2)s$. Then
$e=c\big(1-(1-c)cs\big)\in cR\subseteq aR$, and that $1-e=(1-c)(1+cs)\in
(1-c)R\subseteq bR$. We easily check that $eR+(1-e)R=R, e(1-e)\in J(R)$ and $a|e,b|1-e$. According to Theorem 3.5, $R$ is an elementary divisor ring.\hfill$\Box$

\vskip4mm \hspace{-1.8em} {\bf Example 3.7.}\ \ Let $R:={\Bbb Z}_{(2)}\bigcap {\Bbb Z}_{(3)}=\{ \frac{m}{n}~|~m,n\in {\Bbb Z}, 2,3\nmid n\}$. Then
 $R$ is a feckly clean B\'{e}zout duo-domain, and then it is an elementary divisor domain.

\vskip4mm We say that $c\in R$ is feckly adequate if for any $a\in R$ there exist
some $r,s\in R$ such that $(1)$ $c\equiv rs~(mod~J(R))$; $(2)$ $rR+aR=R$; $(3)$
$s'R+aR\neq R$ for each non-invertible divisor $s'$ of $s$. The following is an generalization of [3, Lemma 3.1] in the duo case.

\vskip4mm \hspace{-1.8em} {\bf Lemma 3.8.}\ \ {\it Every feckly adequate element in B\'{e}zout duo-domains is stable.}
\vskip2mm\hspace{-1.8em} {\it Proof.}\ \ Let $R$ be a B\'{e}zout duo-domain, and let $a\in R$ be feckly adequate. Set $S=R/aR$.
Let $\overline{b}\in S$. Then there exist $r,s\in R$
such that $a\equiv rs (mod~J(R)), rR+bR=R$ and $s'R+bR\neq R$ for any noninvertible
divisor $s'$ of $s$. Hence, $\overline{a}\equiv \overline{rs} \big(mod~J(S)\big)$, i.e., $\overline{rs}\in J(S)$.
Clearly, $\overline{r}S+\overline{b}S=S$. Let
$\overline{t}\in S$ be a noninvertible divisor of
$\overline{s}$. As in the proof of [3, Lemma 3.1], we see that $sR+tR\neq R$. Since $R$ is a B$\acute{e}$zout ring, we have a noninvertible $u\in R$ such that
$sR+tR=uR$. We infer that $u$ is a noninvertible divisor of $s$.
Hence, $uR+bR\neq R$. This proves that
$\overline{u}S+\overline{b}S\neq S$; otherwise,
there exist $x,y,z\in R$ such that $ux+by=1+az$. This implies that
$ux+by=1+w'z+rsz=1+w'z+ucrz$ for $c\in R$ and $w'\in J(R)$, because $a-rs\in J(R)$. Hence, $u(x-crz)(1+w'z)^{-1}+by(1+w'z)^{-1}=1$, a
contradiction. Thus $\overline{t}S+\overline{b}S\neq
S$. Therefore, $\overline{0}\in S$ is feckly adequate.

Let $x\in
S$. by the preceding discussion, we have some $r,s\in R$ such that $rs\in J(S)$,
where $rS+xS=S$ and $s'S+xS\neq S$ for any
noninvertible divisor $s'$ of $s$. Similarly, we have $rS+sS=S$.
Since $rS+xS=S$, we get $rS+sxS=S$, as $S$ is duo.
Write
$rc+sxd=1$ in $S$. Set
$e=rc$. Then $e^2-e=(rc)^2-rc=-(rc)(sxd)\in rsS\subseteq J(S)$. Let $u$ be an arbitrary invertible element of $S$.
As in the proof of [3, Lemma 3.1], we verify that $(u-ex)S+rsS=S$, as $sS=Ss$. Since $rs\in J(S)$, we get $u-ex\in U(S)$. This implies that
$\big(x-(1-e)x\big)-u=ex-u\in U(S)$. Furthermore, $S$ has stable range 1.
We infer that $x-(1-e)x\in J(S)$. Clearly, $1-e=sxd\in xR$. Therefore $\overline{x}=\overline{x(sd)x}$ in $S/J(S)$, and so $S/J(S)$ is regular.

As $R$ is duo, so is $S/J(S)$. Hence, $S$ has stable range 1, by [1, Corollary 1.3.15]. Therefore $a\in R$ is stable, as required.\hfill$\Box$

\vskip4mm A B\'{e}zout domain $R$ is a Zabavsky domain if, $aR+bR=R$ with $a,b\in R\Longrightarrow ~\exists y\in R$ such that $a+by\in R$ is feckly adequate. For instance, every domain of stable range 1 is a Zabavsky domain. We can derive

\vskip4mm \hspace{-1.8em} {\bf Theorem 3.9.}\ \ {\it Every Zabavsky duo-domain is an elementary divisor ring.}\vskip2mm\hspace{-1.8em} {\it Proof.}\ \ Suppose that $aR+bR=R$ with $a,b\in R$. Then there exists a $y\in R$ such that $a+by\in R$ is adequate. In view of Lemma 3.8, $a+by\in R$ is stable. Therefore $R$ is locally stable. This completes the proof, by Theorem 3.2.\hfill$\Box$

\vskip4mm A ring $R$ is called adequate if every nonzero element in $R$ is adequate. We can extend [13, Theorem 1.2.13] as follows:

\vskip4mm \hspace{-1.8em} {\bf Corollary 3.10.}\ \ {\it Every adequate B\'{e}zout duo-domain is an elementary divisor domain.}
\vskip2mm\hspace{-1.8em} {\it Proof.}\ \ Let $R$ be a adequate B\'{e}zout duo-domain. If $aR+bR=R$ with $a,b\in R$, then $a$ or $a+b\neq 0$. Hence, $a+by\in R$ is adequate, where $y=0$ or $1$. This implies that $R$ is a Zabavsky duo-domain. The result follows by Theorem 3.9.\hfill$\Box$

\vskip4mm \hspace{-1.8em} {\bf Theorem 3.11.}\ \ {\it Let $R$ be a B\'{e}zout duo-domain. If for any $a\in R$ either $a$ or $1-a$ is adequate, then $R$ is an elementary divisor domain.}
\vskip2mm\hspace{-1.8em} {\it Proof.}\ \ Given $aR+bR=R$ with $a,b\in R$. Write $ax+by=1$ for some $x,y\in R$. Then
$a(x-y)+(a+b)y=1$. By hypothesis, either $a(1-x)$ or $1-a(1-x)$ is adequate. In view of Lemma 3.8, $a(1-x)$ or $1-a(1-x)$ is stable. If $a(1-x)\in R$ is stable, then $R/a(1-x)R$ has stable range 1. As $R/aR\cong \frac{R/a(1-x)}{aR/a(1-x)R}$, we easily check that $R/aR$ has stable range 1. If $(a+b)y\in R$ is stable, similarly, $R/(a+b)R$ has stable range 1. Therefore $R/(a+by)$ has stable range 1 where $y=0$ or $1$. Hence, $R$ is locally stable, whence the result by Theorem 3.2.\hfill$\Box$\hfill$\Box$

\vskip4mm \hspace{-1.8em} {\bf Corollary 3.12.}\ \ {\it Let $R$ be a B\'{e}zout duo-domain. If for any $a\in R$ either $J(\frac{R}{aR})=0$ or $J(\frac{R}{(1-a)R})=0$, then $R$ is an elementary divisor domain.}
\vskip2mm\hspace{-1.8em} {\it Proof.}\ \ Let $a\in R$. If $a=0$, then $1-a\in R$ is adequate. If $a=1$, then $a\in R$ is adequate. We now assume that $a\neq 0,1$. By hypothesis, $J(\frac{R}{aR})=0$ or $J(\frac{R}{(1-a)R})=0$. In view of [11, Theorem 17], $a$ or $1-a$ is sqaure-free. It follows by [11, Proposition 16] that $a$ or $1-a$ is adequate in $R$. Therefore we complete the proof, by Theorem 3.11.\hfill$\Box$

\vskip4mm A ring $R$ is called homomorphically semiprimitive if every ring homomorphic image (including $R$) of $R$ has zero Jacobson radical. Von Neumann regular rings are clearly homomorphically semiprimitive. But the converse is not true in general, for example, let $R=W_1[F]$ be the first Weyl algebra over a field $F$ of characteristic zero. As an immediate consequence of Corollary 3.12, we have

\vskip4mm \hspace{-1.8em} {\bf Corollary 3.13.}\ \ Let $R$ be a B\'{e}zout duo-domain. If $R$ is homomorphically semiprimitive, then
$R$ is an elementary divisor domain.

\vskip15mm\bc{\bf 4. Rings of Gelfand Range 1}\ec

\vskip4mm A ring $R$ is a PM ring if every prime ideal of $R$ contains in only one maximal ideal. An element $a$ of
a B\'{e}zout duo-domain $R$ is PM if, $\frac{R}{aR}$ is a PM ring. We easily extend [4, Theorem 4.1] to duo-rings.

\vskip4mm \hspace{-1.8em} {\bf Lemma 4.1.}\ \ Let $R$ be
a duo-ring. Then the following are equivalent:
\begin{enumerate}
\item [(1)]{\it $R$ is a PM ring.} \vspace{-.5mm}
\item [(2)]{\it $a+b=1$ in $R$ implies that $(1+ar)(1+bs)=0$ for some $r,s\in R$.}
\end{enumerate}

\vskip4mm \hspace{-1.8em} {\bf Lemma 4.2.}\ \ {\it Let $R$ be a B\'{e}zout duo-domain, and let $a\in R$. Then the following are equivalent:}
\vspace{-.5mm}
\begin{enumerate}
\item [(1)]{\it $a\in R$ is PM.}
\vspace{-.5mm}
\item [(2)]{\it If $bR+cR=R$, then there exist $r,s\in R$ such that $a=rs, rR+bR=sR+cR=R$.}
\vspace{-.5mm}
\item [(3)]{\it If $aR+bR+cR=R$, then there exist $r,s\in R$ such that $a=rs, rR+bR=sR+cR=R$.}
\end{enumerate}
\vspace{-.5mm} {\it Proof.}\ \ $(1)\Rightarrow (3)$ Suppose $aR+bR+cR=R$. Then $\overline{b}(R/aR)+\overline{c}(R/aR)=R/aR$. By virtue of Lemma 4.1, we can find some $p,q\in R$ such that
$\overline{(1+bp)(1+cq)}=\overline{0}$ in $R/aR$. Set $k=1+bp$ and $l=1+cq$. Then $kR+bR=lR+cR=R$ and that $kl\in aR$.

Write $kl=ap$ for some $p\in R$. Since $R$ is a B\'{e}zout duo-domain, we have an $r\in R$ such that $kR+aR=rR$.
Thus, there are some $s,t\in R$ such that $a=rs, k=rt$ and $sR+tR=R$. Clearly, $rR+bR=R$. It follows that $rtl=rsp$, and so $tl=sp$. Write $lu+cv=1$. Then $tlu+ctv=spu+ctv=t$. But $sR+tR=R$, we get $sR+cR=R$, as required.

$(3)\Rightarrow (2)$ This is obvious.

$(2)\Rightarrow (1)$ Suppose $\overline{b+c}=\overline{1}$ in $R/aR$. Then $b+c+ax=1$ for some $x\in R$. By hypothesis,
there exist $r,s\in R$ such that $a=rs, rR+bR=sR+(c+ax)R=R$. Write $rk+bl=sp+(c+ax)q=1$. Then $\overline{(1-bk)(1-cq)}=\overline{rskp}=\overline{0}$ in $R/aR$.
Therefore $R/aR$ is a PM ring, by Lemma 4.1. That is, $a\in R$ is PM, as asserted.\hfill$\Box$

\vskip4mm A B\'{e}zout duo-ring has Gelfand range 1 if, $aR+bR=R$ with $a,b\in R$ implies that there exists a $y\in R$ such that $a+by\in R$ is PM. We come now to the main result of this section:

\vskip4mm \hspace{-1.8em} {\bf Theorem 4.3.}\ \ {\it Every B\'{e}zout duo-domain of Gelfand range 1 is an elementary divisor ring.}\vskip2mm\hspace{-1.8em} {\it Proof.}\ \ Let $R$ be a B\'{e}zout duo-domain of Gelfand range 1. Suppose that $aR+bR+cR=R$ with $a,b,c\in R$. Write $ax+by+cz=1$ with $x,y,z\in R$.
Then we can find some $k\in R$ such that
$w:=b+(ax+cz)k\in R$ is PM.
Hence, $$\left(
\begin{array}{cc}
1&0\\
xk&1
\end{array}
\right)A\left(
\begin{array}{cc}
1&0\\
zk&1
\end{array}
\right)
=\left(
\begin{array}{cc}
a&0\\
w&c
\end{array}
\right).$$
Clearly, $wR+aR+cR=R$. In view of Lemma 4.2, there exist $r,s\in R$ such that $w=rs, rR+aR=sR+cR=R$. Write $sp'+cl=1$. Then $rsp'+rcl=r$. Set $q'=rl$. Then
$wp'+cq'=r$. Write $p'R+q'R=dR$. Then $p'=dp, q'=dq$ and $pR+qR=R$. Hence, $(wp+cq)R+aR=R$. Clearly, $cR+pR=qR+pR=R$, we get $cqR+pR=R$. This implies that
$(wp+cq)R+pR=R$. Thus, we have $apR+(wp+cq)R=R$. Write $aps+(wp+cq)t=1$. Then $p(as+wt)+q(ct)=1$. Hence,
$$\left(
\begin{array}{cc}
s&t\\
wp+cq&-ap
\end{array}
\right)\left(
\begin{array}{cc}
a&0\\
w&c
\end{array}
\right)\left(
\begin{array}{cc}
p&-ct\\
q&as+wt
\end{array}
\right)=\left(
\begin{array}{cc}
1&0\\
0&*
\end{array}
\right),$$ where $\left|
\begin{array}{cc}
s&t\\
wp+cq&-ap
\end{array}
\right|=\left|
\begin{array}{cc}
p&-ct\\
q&as+wt
\end{array}
\right|=1$. As in the proof of Theorem 3.2, we have $p,q\in R$ such that $(pa+qb)R+qcR=R$. This completes the proof, by Lemma 2.7.\hfill$\Box$

\vskip4mm \hspace{-1.8em} {\bf Corollary 4.4 [14, Theorem 3].}\ \ Every commutative B\'{e}zout domain of Gelfand range 1 is an elementary divisor ring.

\vskip4mm We note that Theorem 4.3 can not be extended to any rings with zero divisor as the following shows.

\vskip4mm \hspace{-1.8em} {\bf Example 4.5.}\ \ Let $R^{+}$ be the subset $\lbrace(x,0): x>0\rbrace$ of the plane $S^{+}=\lbrace(x, sin\pi x) : x\geq 0\rbrace$ and let  $X=S^{+}\bigcup R^{+}$. If $R=C(X-\beta(X))$ where $\beta(X)$ the largest compact Hausdorff space generated by $X$ in the sense that any map from $X$ to a compact Hausdorff space factors through $\beta(X)$, and $C(X-\beta(X)$ is the ring of all continuous functions on $X-\beta(X)$ In view of [7], $C(X-\beta(X))$ is a B\'{e}zout duo-ring. On the other hand, $C(X-\beta(x))$ is a PM ring (cf. [4]). Thus $C(X-\beta(x))$ is a B\'{e}zout duo-ring of Gelfand range 1. But it is not Hermite ring (cf. [7]), and therefore $R$ is not an elementary divisor ring.

\vskip4mm \hspace{-1.8em} {\bf Lemma 4.6.}\ \ {\it Let $I$ be an ideal of a PM ring $R$. Then $R/I$ is a PM ring.}
\vskip2mm \hspace{-1.5em}{\it Proof.}\ \ Let $P$ be a prime ideal of $R$. Then $P=Q/I$ where $Q$ is a prime ideal of $R$.
Since $R$ is a PM ring, there exists a unique maximal $M$ of $R$ such that $Q\subseteq M$. Hence,
$P\subseteq M/I$ where $M/I$ is a maximal of $R/I$. The uniqueness is easily checked, and therefore $R/I$ is PM.\hfill$\Box$

\vskip4mm \hspace{-1.8em} {\bf Theorem 4.7.}\ \ {\it Let $R$ be B\'{e}zout duo-domain. Suppose for any $a\in R$, either $a$ or $1-a$ is PM. Then $R$ is an elementary divisor ring.}\vskip2mm \hspace{-1.5em}{\it Proof.}\ \ Given $aR+bR=R$ with $a,b\in R$. Write $ax+by=1$ for some $x,y\in R$. As in the proof of Theorem 3.11, we see that $a(1-x)$ or $1-a(1-x)$ is PM. If $a(1-x)\in R$ is PM, then $R/a(1-x)R$ is PM. Since $R/aR\cong \frac{R/a(1-x)}{aR/a(1-x)R}$, it follows by Lemma 4.6 that $R/aR$ is PM. If $(a+b)y\in R$ is PM, $R/(a+b)R$ is PM. Hence, $R/(a+by)$ is PM, and therefore $R$ is an elementary divisor ring, by Theorem 4.3.\hfill$\Box$

\vskip4mm \hspace{-1.8em} {\bf Corollary 4.8.}\ \ Every B\'{e}zout PM duo-domain is an elementary divisor ring.

\vskip15mm\bc{\bf 5. Countability Conditions}\ec

\vskip4mm Let $R$ be a ring and $a\in R$. We use $mspec(a)$ to denote the set of all
maximal ideals which contain $a$. We now determine elementary divisor domains by means of the countability of the sets of related
elements.

\vskip4mm \hspace{-1.8em} {\bf Lemma 5.1.}\ \ {\it Let $R$ be a B\'{e}zout domain, and let $aR+bR+cR=R$ with $a,b,c\in R$. If $mspec(a)$ is at most countable, then
there exist $p,q\in R$ such that $(pa+qb)R+qcR=R$.}
\vskip2mm\hspace{-1.8em} {\it Proof.}\ \ In view of [13, Theorem 1.2.6], $R$ is an Hermite ring.

Case I. $mspec(a)$ is empty. Then $a\in U(R)$.
Hence, $$\left(
\begin{array}{cc}
a&b\\
0&c
\end{array}
\right)\left(
\begin{array}{cc}
1&-a^{-1}b\\
0&1
\end{array}
\right)=\left(
\begin{array}{cc}
a&0\\
0&c
\end{array}
\right).$$

Case II. $mspec(a)\neq\emptyset$. Set $mspec(a)=\{ M_1,M_2,M_3,\cdots \}$. If $b\in
M_1$, then $c\not\in M_1$; hence, $b+c\not\in M_1$. By an
elementary transformation of column, $b$ can be replaced by the
element $b+c$. Thus, we may assume that $b\not\in M_1$. As $R$ is a Hermite ring, there exists $P_1\in GL_2(R)$ such that
$P_1A=\left(
\begin{array}{cc}
a_1&b_1\\
0&c_1
\end{array}
\right)$. This implies that $aR+bR=a_1R$. If $a_1\in U(R)$, then
$$P_1A\left(
\begin{array}{cc}
1&-a_1^{-1}b_1\\
0&1
\end{array}
\right)=diag(a_1,c_1)$$ is a diagonal matrix.

Assume that $a_1\not\in U(R)$. Then $mspec(a_1)\neq \emptyset$. As
$mspec(a_1)\subseteq mspec(a)$, we may assume that $a_1\in M_2$. If
$b_1\in M_2$, then $c_1\not\in M_2$ since $a_1R+b_1R+c_1R=R$. Hence,
$b_1+c_1\not\in M_2$. By an elementary transformation of row, $b_1$
can be replaced by the element $b_1+c_1$. Thus, we may assume that
$b_1\not\in M_2$. Since $R$ is a Hermite ring, there exists $Q_1\in
GL_2(R)$ such that $P_1AQ_1=\left(
\begin{array}{cc}
a_2&0\\
b_2&c_2
\end{array}
\right)$. Clearly, $a_1R+b_1R=a_2R$. Moreover, we see that
$a\in M_1,a_1\in M_2\backslash M_1, a_2\in M_3\backslash
\big(M_1\bigcup M_2\big)$. By iteration of this process, there is a
collection of matrices of the form
$$P_kAQ_k=\left(
\begin{array}{cc}
a_{i}&*\\
0&*
\end{array}
\right).$$ If there exists some $a_{i}\in U(R)$, then $A$ admits a
diagonal reduction. Otherwise, we
get an infinite chain of ideals $$aR\subsetneqq a_1R\subsetneqq
a_2R\subsetneqq \cdots .$$ Further, each $a_n\not\in
M_{n}(n\in {\Bbb N}$. Let $I=\bigcup\limits_{i=1}^{\infty}a_iR$. Then
$I\neq R$; hence, there exists a maximal ideal $J$ of $R$ such
that $I\subseteq J\subsetneq R$. This implies that $a\in J$, and so
$J\in mspec(a)$. Thus, $J=M_k$ for some $k\in \mathbb{N}$. As
$a_k\not\in M_k$, we deduce that $a_k\not\in J$, a contradiction.
Therefore $A$ admits a diagonal reduction. As in the proof of Theorem 3.2, there exist $p,q\in R$ such that $(pa+qb)R+qcR=R$, as required.\hfill$\Box$

\vskip4mm We say that an element $a$ of a ring $R$ is a non-P element if, $\frac{R}{aR}$ do not satisfy the property $P$, where $P$ is some ring theoretical property. We now derive

\vskip4mm \hspace{-1.8em} {\bf Theorem 5.2.}\ \ {\it Let $R$ be a B\'{e}zout duo-domain. Then $R$ is an elementary divisor ring if}
\vspace{-.5mm}
\begin{enumerate}
\item [(1)]{\it the set of nonstable elements in $R$ is at most countable, or}
\vspace{-.5mm}
\item [(2)]{\it the set of non-PM elements in $R$ is at most countable.}
\end{enumerate}
\vspace{-.5mm} {\it Proof.}\ \ $(1)$ If $R$ has stable range 1, then it is an elementary divisor ring, by Lemma 3.4. We now assume that $R$ has no stable range 1. Suppose that $aR+bR+cR=R$ with $a,b,c\in R$. If $a\in R$ is stable, as in the proof of Theorem 3.2, we can find some $p,q\in R$ such that $(pa+qb)R+qcR=R$. If
$a\in R$ is nonstable, then $a\neq 0$; otherwise, $R$ has stable range 1. Since $R$ is a domain, we see that
$a^m\neq a^n$ for any distinct $m,n\in {\Bbb N}$. If $mspec(a)$ is uncountable, then $\{ a,a^2,\cdots ,a^k,\cdots \}$ is uncountable as a subset of $mspec(a)$. This gives a contradiction. This implies that $mspec(a)$ is at most countable.
By virtue of Lemma 5.1, there exist $p,q\in R$ such that $(pa+qb)R+qcR=R$.
Therefore proving $(1)$, in terms of Lemma 2.7.

$(2)$ Suppose that $aR+bR+cR=R$ with $a,b,c\in R$.
If $a\in R$ is stable, similarly to Theorem 3.2, we have some $p,q\in R$ such that $(pa+qb)R+qcR=R$. If
$a\in R$ is a non-PM element, by hypothesis, $mspec(a)$ is at most countable. In light of Lemma 5.1, we can find $p,q\in R$ such that $(pa+qb)R+qcR=R$.
Thus, we obtain the result by Lemma 2.7.\hfill$\Box$

\vskip4mm Recall that an ideal of a duo ring $R$ is a maximally nonprincipal ideal if it is maximal in the set of nonprincipal ideals of $R$ with respect to the ordering by set inclusion.

\vskip4mm \hspace{-1.8em} {\bf Lemma 5.3.}\ \ {\it Let $R$ be a B\'{e}zout domain, and let $aR+bR+cR=R$ with $a,b,c\in R$. If $a$ does not belong to any maximally nonprincipal ideal, then
there exist $p,q\in R$ such that $(pa+qb)R+qcR=R$.}
\vskip2mm\hspace{-1.8em} {\it Proof.}\ \ Set $A=\left(
\begin{array}{cc}
a&0\\
b&c
\end{array}
\right)\in M_2(R)$. Since $a$ does not belong to any maximally nonprincipal ideal, it follows by [7, Theorem 1] that $A$ admits a diagonal reduction. Therefore we complete the proof, as in the proof of Theorem 3.2.\hfill$\Box$

\vskip4mm \hspace{-1.8em} {\bf Theorem 5.4.}\ \ {\it Let $R$ be a B\'{e}zout duo-domain. Then $R$ is an elementary divisor ring if}
\vspace{-.5mm}
\begin{enumerate}
\item [(1)]{\it every nonstable element in $R$ is contained in a maximal nonprincipal ideal, or}
\vspace{-.5mm}
\item [(2)]{\it every non-PM element in $R$ is contained in a maximal nonprincipal ideal.}
\end{enumerate}
\vspace{-.5mm} {\it Proof.}\ \ $(1)$ Suppose that $aR+bR+cR=R$ with $a,b,c\in R$. If $a\in R$ is stable, similarly to Theorem 3.2, there are some $p,q\in R$ such that $(pa+qb)R+qcR=R$. If
$a\in R$ is nonstable, by hypothesis, $a\in R$ is contained in a maximal nonprincipal ideal. In light of Lemma 5.3,
$(pa+qb)R+qcR=R$ for some $p,q\in R$. It follows by Lemma 2.7 that $R$ is an elementary divisor ring.

$(2)$ Analogously to Theorem 4.3 and the preceding discussion, we complete the proof, by Lemma 5.3 and Lemma 2.7.\hfill$\Box$

\vskip4mm \hspace{-1.8em} {\bf Corollary 5.5.}\ \ {\it Let $R$ be a B\'{e}zout duo-domain. Then $R$ is an elementary divisor ring if}
\vspace{-.5mm}
\begin{enumerate}
\item [(1)]{\it every nonlocal element in $R$ is contained in a maximal nonprincipal ideal, or}
\vspace{-.5mm}
\item [(2)]{\it every nonregular element in $R$ is contained in a maximal nonprincipal ideal.}
\end{enumerate}
\vspace{-.5mm} {\it Proof.}\ \ This is obvious by Theorem 5.4, as every local duo-ring and every regular duo-ring are both PM rings. \hfill$\Box$

\vskip15mm \bc{\large\bf References}\ec \vskip4mm {\small \re{1} H.
Chen, {\it Rings Related Stable Range Conditions}, Series in
Algebra 11, World Scientific, Hackensack, NJ, 2011.

\re{2} H. Chen; H. Kose and Y. Kurtulmaz, On feckly clean rings,
{\it J. Algebra Appl.}, {\bf 14}(2015), 1550046(15 pages). DOI:10.1142/S0219498815500462.

\re{3} H. Chen and M. Sheibani, Elementary matrix reduction over Zabavsky rings,
{\it Bull. Korean Math. Soc.}, to appear.

\re{4} M. Contessa, On pm-rings, {\it
Comm. Algebra}, {\bf 10}(1982), 93--108.

\re{5} O.V. Domsha and I.S. Vasiunyk, Combining local and adequate rings, {\it Book of abstracts of the International Algebraic Conference}, Taras Shevchenko National University of Kyiv, Kyiv, Ukraine, 2014.

\re{6} A. Gatalevich and B.V. Zabavsky, Noncommutative elementary divisor domains,
{\it J. Math. Sci.}, {\bf 96}(1999),
3013--3016.

\re{7} L. Gillman and M. Henriksen, Rings of continuous functions in which every finitely generated ideal is principal,
{\it Trans. Amer. Math. Soc.}, {\bf 82}(1956), 362--365.

\re{8} T.Y. Lam and A.S. Dugas, Quasi-duo rings and stable range descent,
{\it J. Pure Appl. Algebra}, {\bf 195}(2005), 243--259.

\re{9} W.W. McGovern, B\'{e}zout rings with almost stable range
$1$, {\it J. Pure Appl. Algebra}, {\bf 212}(2008), 340--348.

\re{10} M. Roitman, The Kaplansky condition and rings of almost
stable range $1$, {\it Proc. Amer. Math. Soc.}, {\bf 141}(2013),
3013--3018.

\re{11} O.S. Sorokin, Finite homomorphic images of B\'{e}zout duo-domain,
{\it Carpathian Math. Publ.}, {\bf 6}(2014), 360--366.

\re{12} S.H. Sun, Noncommutative rings
in which every prime ideal is contained in a unique maximal ideal,
{\it J. Pure Appl. Algebra}, {\bf 76}(1991), 179--192.

\re{13} B.V. Zabavsky, Diagonal Reduction of Matrices over Rings,
Mathematical Studies Monograph Series, Vol. XVI, VNTL Publisher, 2012.

\re{14} B.V. Zabavsky, Conditions for stable range of an elementary divisor rings, arXiv:1508.07418v1 [math.RA] 29 Aug 2015.

\re{15} B.V. Zabavsky and N.Y. Komarnitskii, Distributive elementary divisor domains, {\it Ukrainskii Math. Zhurnal}, {\bf 42}(1990),
1000--1004.

\end{document}